\documentclass[10pt]{article}
\usepackage{amsfonts,amsmath,amsthm,amssymb}
\usepackage{graphics, epsfig}
\usepackage{color}
\usepackage{appendix}
\usepackage{verbatim}
\usepackage{ulem}
\usepackage[makeroom]{cancel}
\usepackage[margin=1.5in]{geometry}

 \usepackage[usenames,dvipsnames]{pstricks}
 \usepackage{pst-grad} 
 \usepackage{pst-plot} 
\allowdisplaybreaks

\newtheorem{proposition}{Proposition}[section]
\newtheorem{theorem}{Theorem}[section]
\newtheorem{lemma}{Lemma}[section]

\newtheorem{remark}{Remark}[section]

\numberwithin{equation}{section}
\numberwithin{theorem}{section}
\numberwithin{proposition}{section}
\numberwithin{lemma}{section}
\numberwithin{remark}{section}
\newcommand{\gm}{\gamma}

\newcommand{\z}{\zeta}

\newcommand{\essosc}{\operatornamewithlimits{ess\,osc}}

\newcommand{\dist}{\operatorname{dist}}

\newcommand{\intl}{\int\limits}

\def\Xint#1{\mathchoice
    {\XXint\displaystyle\textstyle{#1}}%
    {\XXint\textstyle\scriptstyle{#1}}%
    {\XXint\scriptstyle\scriptscriptstyle{#1}}%
    {\XXint\scriptscriptstyle\scriptscriptstyle{#1}}%
    \!\int}
\def\XXint#1#2#3{\setbox0=\hbox{$#1{#2#3}{\int}$}
    \vcenter{\hbox{$#2#3$}}\kern-0.5\wd0}
\def\bint{\Xint-}
\def\dashint{\Xint{\raise4pt\hbox to7pt{\hrulefill}}}
\def\dashiint{\bint\kern-0.15cm\bint}

\newcommand{\ovl}[3]{\int_{#1}^{#2}\kern-#3pt\raise4pt\hbox to7pt{\hrulefill}\ }

\newcommand{\ovll}[3]{\intl_{#1}^{#2}\kern-#3pt\raise4pt\hbox to7pt{\hrulefill}\ }

\newcommand{\tvl}[2]{\iint_{#1}\kern-#2pt\raise4pt\hbox to7pt{\hrulefill}\ }

\newcommand{\bye}{\end{document}}

\newcommand{\tvls}[2]{\iint_{#1}\kern-#2pt\raise4pt\hbox to15pt{\hrulefill}\ }

\begin{document}
\title{A new short proof of regularity for local weak solutions for a certain class of singular parabolic equations}
\author{
\\
\\
\it{Simone Ciani \& Vincenzo} Vespri \\ \\ Universit\`a degli Studi di Firenze,\\ Dipartimento di Matematica e Informatica "Ulisse Dini"\\   $simone.ciani@unifi.it$ \& $vincenzo.vespri@unifi.it$
}
\date{}
\maketitle
\vskip.4truecm
\begin{abstract} \noindent
We shall establish the interior H\"older continuity for locally bounded weak
solutions to a class of parabolic singular equations whose prototypes are
\begin{equation} \label{p-laplacean}
u_t= \nabla \cdot \bigg( |\nabla u|^{p-2} \nabla u \bigg), \quad \text{ for } \quad 1<p<2,
\end{equation}  and
\begin{equation}\label{doublysingular}
 u_{t}- \nabla \cdot ( u^{m-1} | \nabla u |^{p-2} \nabla u ) =0 , \quad \text{for} \quad m+p>3-\frac{p}{N},
\end{equation}
via a new and simplified proof using recent techniques on expansion of positivity and $L^{1}$-Harnack estimates.
\vskip.2truecm
\noindent
{\bf{MSC 2020:}}
35K67, 35K92, 35B65
\vskip.2truecm
\noindent
{\bf{Key Words}}: 
Singular Parabolic Equations, $p$-Laplacean, Doubly Nonlinear, H\"older Continuity, Intrinsic Scaling, Expansion of Positivity.
\\
\begin{flushright}
\it{To celebrate Umberto Mosco's 80th genethliac}
\end{flushright}
\end{abstract}

\section{Introduction}\label{Intro}
Equations of the kind of \eqref{p-laplacean} are termed singular since, the modulus of ellipticity 
$|\nabla u|^{p-2}$ becomes infinitely big as the weak gradient of the function $|\nabla u|$ approaches zero. Regularity theory and in particular the study of H\"older continuity for such singular parabolic equations has been pioneered by Y.Z. Chen and E. Di Benedetto in \cite{Chen-DiBe88}, \cite{Chen-DiBe92}. The singular approach is more difficult than the degenerate one, i.e. when $p>2$ where the modulus of ellipticity tends to vanish. A detailed study for the class of degenerate parabolic equations of $p$-Laplacean type has been extensively treated in the monograph \cite{DB}. The method developed to achieve the continuity of local weak solutions of both degenerate and singular equations of these kind bears the name of intrinsic scaling. This approach was introduced by E. DiBenedetto (see the monograph \cite{DB}, see also \cite{Urbano}) and its name comes from the fact that the diffusion processes in the equations evolve in a time
scale determined instant by instant by the solution itself, so that, loosely speaking, they can be regarded as the heat equation in their own intrinsic time configuration. To overcome the difficulties of this approach, it was introduced a more geometrical method named expansion of positivity. It was initially developed in the degenerate case for the study of Harnack inequality (see \cite{DBGV-acta}) and then used to give a more direct proof of regularity in \cite{Surnachev}. In the singular case, the expansion of positivity was proved in \cite{Proceedings}, and it was simply used in \cite{DBGV-mono} to avoid the use of a very technical Lemma that is central in the proof of \cite{Chen-DiBe92}.
The aim of this paper is to use the full potentiality of the expansion of positivity Lemma in order to give a more direct and geometrical proof of regularity of solutions to singular equations of the kind of \eqref{p-laplacean}.

 $$$$
The method we present here can also be implemented for solutions to doubly nonlinear equations of the kind \eqref{doublysingular}. The expansion of positivity was proved in \cite{Vespri-Sosio-Fornaro2} (see also \cite{Vespri-Sosio-Fornaro1} and \cite{Matias}). Equations as \eqref{doublysingular} are the natural bridge between the porous media equations and $p$-Laplace type ones. They constituted and still constitute an hard challenge from the mathematical point of view, because many questions (also of regularity) are still open. The term doubly nonlinear refers to the fact that the diffusion part depends nonlinearly both
on the gradient and the solution itself. These equations have been introduced by J.L.Lions in \cite{Lions} and they describe several physical phenomena; see the survey of A.S.Kalashnikov \cite{Kalashnikov} for more details, see also the following papers for a non-comprehensive surveys on this argument, \cite{Kuusi1}, \cite{Kuusi2}, \cite{Matias-Singer} and \cite{Sturm}.
In this paper we take as a starting point the recent extensive study made in \cite{Matias} and we also refer to it for a self-contained introduction to the regularity theory for doubly nonlinear equations.$$$$
Let us sketch the strategy for the proof of H\"older continuity in the case of doubly nonlinear equations; the $p$-Laplacean case is easier. Let us recall that we follow the De Giorgi's approach (\cite{DG}) where the H\"older continuity was proved via the reduction of oscillation.\\
If $u$ is the solution, for sake of simplicity, assume that the solution $u$ satisfies $0 \leq u \leq 1$. Let $Q$ be a cylinder, and we state an alternative on the measure of the set where the solution $u$ is greater than $\frac{1}{2}$. Either the measure of this set is greater than a sizeable portion of the cylinder itself or this measure is smaller. We have two alternatives.\\
Assume that $[u>\frac{1}{2}] \cap Q|\leq \nu |Q|$, where $\nu$ is a suitable constant in $(0,1)$ to be chosen. For $\nu$ small enough, it is possible to apply a De Giorgi's result (the so-called Critical Mass Lemma) to get that in a smaller cylinder the solution is smaller than $\frac{3}{4}$, and this implies the reduction of oscillation. \\
If the other alternative happens, i.e. $[u>\frac{1}{2}] \cap Q|> \nu |Q|$, we have that the measure of the set where $u$ is "big" is itself big. Then there is a time level $\bar{t}$ where in the ball $B$ we have $|u(\cdot, \bar{t})> \frac{1}{2}] \cap B| > \nu |B|$. Let us apply an integral Harnack estimate introduced for the first time for the $p$-Laplacean in \cite{Chen-DiBe88} (see also \cite{Annali}) and for the doubly nonlinear case in \cite{Vespri-Sosio-Fornaro1} (see also \cite{Matias}). Thanks to this inequality, the measure information can be extended to any time level in $Q$. Hence we are under the assumptions where we can apply the expansion of positivity Lemma, and so we are able to find a subcylinder $Q' \subset Q$ where the solution is greater than a small constant. In this way, we have a reduction of the oscillation of $u$ and thus the H\"older continuity of the solution is proved.
$$$$
The present paper is organised as follows. In \textsection \ref{S:2} we introduce notations and main results for both the class of equations. In \textsection \ref{S:3} we prove H\"older continuity for local weak solutions to equations of the kind of \eqref{p-laplacean}, 
and finally we devolve \textsection \ref{S:4} to the proof of H\"older continuity for local weak solutions to doubly singular equations as \eqref{doublysingular}.

\section{Notation and Main Results} \label{S:2} $$$$
\subsection{The case of $p$-Laplacean equations}
Let $\Omega$ be an open set in $\mathbb{R}^N$ and for $T>0$ let $\Omega_T$ denote the cylindrical domain $\Omega \times (0,T]$. We denote by $|E|$ the Lebesgue measure of the set $E \subset \mathbb{R}^{N}$ and for a $k \in \mathbb{R}$ by $ [u>k] \cap E$ the set of points of $E$ in which the inequality $u>k$ holds. We write $\nabla u$ the gradient of $u$ taken with respect to the spatial variables, and with $\nabla \cdot {\bf{v}}$ the spatial divergence of a vector field ${\bf{v}}$.
Consider quasi-linear, parabolic differential equations of the form 
\begin{equation} \label{Eq:1:1}
\begin{aligned}
&u \in C_{loc} (0,T; L^{2}_{loc}(\Omega)) \cap L^{p}_{loc}(0,T; W^{1,p}_{loc}(\Omega)),\\
&u_t - \nabla \cdot \bigg( {\bf{A}}(x,t,u,\nabla u) \bigg)= 0  \quad \text{in } \quad D'(\Omega_T)
\end{aligned}\quad \quad \quad  \text{for} \quad 1<p<2
\end{equation} 
where for $C>0$
\begin{equation}
\begin{cases}
{\bf{A}}(x,t,u,\nabla u) \cdot  \nabla u \ge C | \nabla u|^p\\
|{\bf{A}}(x,t,u,\nabla u)| \leq C |\nabla u|^{p-1}
\end{cases}
\end{equation} \noindent 
A measurable function $u$ is a local weak solution of \eqref{Eq:1:1} in $\Omega_T$ if
\begin{equation}
u \in C_{loc} (0,T;L^{2}_{loc}(\Omega)) \cap L^p_{loc}(0,T;W^{1,p}_{loc}(\Omega)),
\end{equation} \noindent
and for every compact subset $K \subset \subset \Omega$ and for every sub-interval $[t_1,t_2] \subset (0,T]$
\begin{equation} \label{solutiondef}
\begin{aligned}
\int_K u \varphi \, dx \bigg{|}_{t_1}^{t_2} +& \int_{t_1}^{t_2} \int_K \{ -u \varphi_t+ |\nabla u|^{p-2} \nabla u \cdot \nabla \varphi \}\, dx d\tau
 = 0 
\end{aligned}
\end{equation} \noindent for all locally bounded testing functions
\begin{equation}
\varphi \in W^{1,2}_{loc} (0,T;L^2_{loc}(K)) \cap L^p_{loc}(0,T; W^{1,p}_{o} (K)), 
\end{equation} \noindent $$$$
For $\rho >0$ let $B_{\rho}$ be the ball of center the origin in $\mathbb{R}^N$ and radius $\rho$ whose boundary is denoted by $\Gamma$. For $y \in \mathbb{R}^N$ let $B_{\rho}(y)$ be the homothetic ball centered at $y$. Let $w_N$ be the measure of the $N$-dimensional unitary ball. Finally for $\rho, l>0$ denote by $Q(l,\rho)= B_{\rho} \times (-l,0]$ the standard cylinder. 

\begin{proposition}[$p$-Laplacean Expansion in Positivity \cite{Proceedings}] \label{P:1:1}
Let $u$ be a non-negative, local, weak solution to \eqref{Eq:1:1}, satisfying
\begin{equation} \label{Eq:1:2}
|[u(\cdot,t)>M] \cap B_{\rho}(y)|>\alpha |B_{\rho}|
\end{equation} \noindent for all times
\begin{equation}
s-\epsilon M^{2-p} \rho^p \leq t \leq s
\end{equation} \noindent for some $M>0$, and $\alpha, \epsilon  \in (0,1)$, and assume that for a fixed number $m \in \mathbb{N}$ it holds
\begin{equation}
B_{8m\rho}(y) \times [s-\epsilon M^{2-p} \rho^p,s] \subset \Omega_T.
\end{equation} \noindent
Then there exist $\sigma \in (0,1)$ and $\epsilon^* \in (0, \frac{1}{2}\epsilon]$ ,wich can be determined a priori, quantitatively only in terms of the data, and the numbers $\alpha, \epsilon, m$, and independent of $M$, such that
\begin{equation}
u(x,t) \ge \sigma M, \quad \text{for all} \quad x \in B_{m \rho}(y),
\end{equation} \noindent for all times
\begin{equation}
s-\epsilon^* M^{2-p} \rho^p <t \leq s.
\end{equation}
\end{proposition}
In addition, for the proof of H\"older continuity, we will need the following estimate from [\cite{DB}, Prop.4.1 pg 193]. The Proposition can be regarded as a weak integral form of a Harnack estimate. That is, the $L^1$-norm of $u(\cdot,t)$ over a ball controls the $L^1$-norm of $u(\cdot,\tau)$ over a smaller ball, for any previous or later time in a suitable interval.
\begin{proposition}[Integral Harnack inequality \cite{DB}] \label{integralHarnack}
Let $u$ be a non-negative weak solution of \eqref{Eq:1:1} and let $1<p<2$. There exists a constant $\gm=\gm(N,p)$ such that
$$\forall(x_0,t_0) \in \Omega_T, \quad \forall \rho>0, \quad \text{such that} \quad B_{4\rho}(x_0) \subset \Omega, \quad \forall t>t_0$$
\begin{equation}
\sup_{t_0\leq \tau \leq t} \int_{B_{\rho}(x_0)} u(x,\tau) dx \leq \gm \inf_{t_0 \leq \tau \leq t} \int_{B_{2\rho} (x_0)} u(x,\tau) dx+ \gm \bigg( \frac{t-t_0}{\rho^{N(p-2)+p}}  \bigg)^{\frac{1}{2-p}}
\end{equation}
\end{proposition}
\begin{remark}
The proof shows that the constant $\gm(N,p)$ deteriorates as $p \rightarrow 2$.
\end{remark}
Finally we state the main theorem as our result.
\begin{theorem} \label{T:1:1}
Let $u$ be a bounded local weak solution of \eqref{Eq:1:1}. Then $u$ is locally H\"older continuous in $\Omega_T$, and there exist constants $\gm >1$ and $\alpha \in (0,1)$ depending only upon the data, such that $\forall K \subset \Omega_T$ compact set,
\begin{equation}
|u(x_1,t_1)-u(x_2,t_2)| \leq \gm ||u||_{\infty, \Omega_T} \bigg(\frac{||u||_{\infty, \Omega_t}^{\frac{2-p}{p}}|x_1-x_2|+|t_1-t_2|^{\frac{1}{p}}}{p-\dist (K;\Gamma)}   \bigg)^{\alpha}
\end{equation} \noindent where $p-\dist$ denotes the intrinsic parabolic distance from $K$ to the parabolic boundary $\Omega_T$, i.e.
\begin{equation}
p-\dist(K;\Gamma):= \inf_{(x,t) \in K, \quad (y,s) \in \Gamma} \bigg( ||u||_{\infty, \Omega_T}^{\frac{2-p}{p}}|x-y|+|t-s|^{\frac{1}{p}}   \bigg)
\end{equation}
\end{theorem} \noindent 
The Theorem \ref{T:1:1} will be proved if reduction of oscillation can be achieved. For sake of completeness we give the explanation to this fact by next Proposition which can be found in \cite{DB} pages 80-81.
\begin{proposition} \label{P:iteration}
Suppose that there exist constants $a, \epsilon^* \in (0,1)$ and $b, \Gamma>1$ that can be determined only in terms of the data, satisfying the following. Construct the sequences
\[ \rho_n= b^{-n} \rho, \quad \quad \rho_0=\rho \quad \quad \quad \forall n= 0,1,2,..
\]
\[
\omega_{n+1}= \max \{ a\omega, \Gamma \rho_n^{\epsilon^*}\}, \quad \quad \omega_0= \omega, \quad \quad \quad \forall n=0,1,2,..
\] and the cylinders
\[
Q_n= Q(\rho_n^p, c_n \rho_n), \quad \text{with} \quad c_n= \omega_n^{\frac{p-2}{p}}, \quad  \quad \quad \forall n=0,1,2,..
\] such that, for all $n=0,1,2,..$ it holds
\[
Q_{n+1} \subset Q_n, \quad \text{and} \quad \essosc_{Q_n} u \leq \omega_n
\]
Then there exist constants $\gm >1$ and $\alpha \in (0,1)$ that can be determined a priori only in terms of the data, such that for all cylinders 
\[
0< r \leq \rho, \quad \quad Q(r^p,c_0 r), \quad \quad c_0= \omega^{\frac{p-2}{p}}
\] holds
\[
\essosc_{Q(r^p,c_o r)} u \leq \gm ( \omega + \rho^{\epsilon^*}) \bigg(\frac{r}{\rho}  \bigg)^{\alpha}
\]
\end{proposition} \noindent
H\"older continuity over compact subsets of $\Omega_T$ is therefore implied by this estimate by a standard covering argument.
$$$$

\subsection{The doubly nonlinear case}
Let us consider the weak solutions to doubly nonlinear equations whose model case is
\begin{equation} \label{doublynonlinear}
    u_{t}- \nabla \cdot ( u^{m-1} | \nabla u |^{p-2} \nabla u ) =0 \quad \quad \text{in}  \quad \Omega_T:= \Omega \times (0,T),
\end{equation} \noindent where $\Omega $ is an open bounded subset of $\mathbb{R}^N$ and 
\begin{equation} \label{range}
p \in (1,2), \quad \quad m>1, \quad \quad \text{and} \quad \quad 2<m+p<3.
\end{equation} \noindent We recall that when $m+p>3$ we are in the degenerate case, while when $m+p=3$ the equation behaves like the heat equation and is called Trudinger's equation, named in this way because introduced by Trudinger in \cite{Trudy}. When $2<m+p<3$ we are in the singular case. The case $m+p=2$ was considered in \cite{Fornaro-Henriquez-Vespri}, where the doubly nonlinear equation has a logarithmic behavior.\\
We will prove the H\"older continuity of solutions for the so-called supercritical range, where we have a consolidated theory developed, i.e.
\begin{equation} \label{supercritical}
3-\frac{p}{N}<m+p<3
\end{equation} \noindent 
The classical theory (see for instance \cite{Matias}) shows that the equation \eqref{doublynonlinear} can be transformed into 
\begin{equation*} \label{formallydoublynonlinear}
u_t- \nabla \cdot (\beta^{1-p} | \nabla u^\beta |^{p-2} \nabla u^{\beta})=0,\quad \quad \quad \beta= \frac{p+m-2}{p-1}
\end{equation*} \noindent This transformation is useful in order to avoid proofs involving the weak gradient of $u$, which has been shown in \cite{Ivanov2} to be existing. The technique that we are applying works perfectly for more general equations of this kind, as
\begin{equation} \label{generalformallydoublynonlinear}
u_t- \nabla \cdot {\bf{A}}(x,t,u,\nabla u^{\beta})=0,
\end{equation} \noindent where ${\bf{A}}$ is a Caratheodory vector field satisfying the conditions
\begin{equation*} \label{Agrowth}
|{\bf{A}}(x,t,s,\zeta)| \leq C_1 | \zeta|^{p-1}$$$$
{\bf{A}}(x,t,s,\zeta) \cdot \zeta \ge C_0 |\zeta|^p
\end{equation*} \noindent
A weak solution for the equation \eqref{generalformallydoublynonlinear} is a non negative function $u:\Omega_T \rightarrow \mathbb{R}$, $u^{\beta} \in L^p(o,T; W^{1,p}(\Omega))$, $u \in L^{\beta+1} ( \Omega_T)$, such that
\begin{equation} \label{definitiondnnl}
\int \int_{\Omega_T} [ {\bf{A}}(x,t,u \nabla u^{\beta}) \cdot \nabla \phi - u \phi_t  ] dxdt =0, \quad \quad \forall \phi \in C^{\infty}_o(\Omega_T)
\end{equation} \noindent 
The proof in this case is different from the $p$-Laplacean one, because $(1-u)$ is not anymore a solution to the previous equation.
Next we give the geometrical setting we use to state the main Lemmata. Let  $0<M<\infty$ and
\begin{equation} 
\theta= (2M)^{3-m-p}
\end{equation} \noindent 
Pick $(\bar{x},\bar{t}) \in \Omega_T$ and suppose that for $\theta$ as above, and a sufficiently small $ 0<\rho<1$ the cylinder
$Q_{\rho}^- (\theta)= (\bar{x},\bar{t})+ B_{\rho} \times (-\theta \rho^p,0]$ is contained in $\Omega_T$.\\ 
The following is a doubly nonlinear version of the Critical Mass Lemma, which is a slight modification of \cite{Matias}.
\begin{lemma}[Critical Mass Lemma \cite{Matias}] \label{DG}
Suppose that $u$ is a weak solution to the equation \eqref{generalformallydoublynonlinear}, and suppose that there exists $M>0$ and $\theta$ defined as above such that the cylinder $Q_{\rho}^-(\theta) \subset \Omega_T$, and it is satisfied
\begin{equation} \label{M}
    \sup_{Q_{\rho}(\theta)} u \leq 2M.
\end{equation} \noindent Then there exists a constant $\nu \in (0,1)$ depending only on the data such that if
\begin{equation} \label{DG-condition}
|[u>M] \cap Q_{\rho}^{-}(\theta)| \leq  \frac{\nu}{(\theta M^{m+p-3})}  | Q_{\rho}(\theta)^- | 
\end{equation} \noindent then we have
$$u \leq \sqrt[\beta]{\frac{3}{2}} M \quad \quad \text{in} \quad Q_{\frac{\rho}{2}}^- (\theta)$$
\end{lemma} \noindent 
Next we state a Lemma of expansion of positivity.
\begin{lemma}[Expansion of Positivity \cite{Matias}]
Suppose that $(x_o,s) \in \Omega_T$ and $u$ is a weak solution of \eqref{generalformallydoublynonlinear}, satisfying for $M>0, \alpha \in (0,1)$
\begin{equation}
|B_{\rho} (x_o) \cap [u(\cdot,s) \ge M]| \ge \alpha |B_{\rho}(x_o)|
\end{equation}\noindent Then there exist $\epsilon,\delta,\eta \in (0,1)$ depending only on the data such that if $$B_{16\rho}(x_o) \times (s,s+\delta M^{3-m-p}\rho^p) \subset \Omega_T$$ then
$$u \ge \eta M \quad \quad \text{in} \quad \text B_{2\rho}(x_o) \times (s+(1-\epsilon) \delta M^{3-m-p} \rho^p,s+\delta M^{3-m-p}\rho^p)$$
\end{lemma} \noindent Finally we recall the following $L^1$-Harnack inequality which was first demonstrated in \cite{Vespri-Sosio-Fornaro1}
\begin{lemma}[Integral Harnack Inequality \cite{Matias}] \label{Harnack DBNL}
Let $u$ be a weak solution of equation \eqref{generalformallydoublynonlinear}. Then there exists $\gm>0$ depending only upon the data such that for all chosen cylinder $B_{2\rho}(y) \times [s,T] \subset \subset \Omega_T$
\begin{equation}
    \sup_{s \leq \tau \leq T} \int_{B_{\rho}(y)} u(x,\tau) dx \leq \gm \bigg{\{} \inf_{s \leq \tau \leq T} \int_{B_{2 \rho}}u(x,\tau) dx + \bigg[ \frac{(T-s)}{\rho^p}\bigg]  ^{\frac{1}{3-m-p}} \rho^N  \bigg{\}}
\end{equation}
\end{lemma} \noindent
Through the previous results, we are able to prove the following theorem.

\begin{theorem} \label{HCDBNL}
Let $u$ be a local weak solution of the equation \eqref{generalformallydoublynonlinear} and let $m,p$ be in the supercritical range \eqref{supercritical}. Then $u$ is locally H\"older continuous in $\Omega_T$, i.e. there exist a H\"older exponent $\alpha \in (0,1)$ depends only on $m,N,p,C_0,C_1$, a constant $\gm>1$, such that $\forall K \subset \Omega_T$ compact set,
\begin{equation}
|u(x_1,t_1)-u(x_2,t_2)| \leq \gm ||u||_{\infty, \Omega_T} \bigg(\frac{||u||_{\infty, \Omega_t}^{\frac{3-m-p}{p}}|x_1-x_2|+|t_1-t_2|^{\frac{1}{p}}}{(m,p)-\dist (K;\Gamma)}   \bigg)^{\alpha}
\end{equation} \noindent where $(m,p)-\dist$ denotes the intrinsic parabolic weighted distance from $K$ to the parabolic boundary $\Omega_T$, i.e.
\begin{equation}
(m,p)-\dist(K;\Gamma):= \inf_{(x,t) \in K, \quad (y,s) \in \Gamma} \bigg( ||u||_{\infty, \Omega_T}^{\frac{3-m-p}{p}}|x-y|+|t-s|^{\frac{1}{p}}   \bigg)
\end{equation}
\end{theorem}
\begin{remark}
We recall that if $u$ is a local weak solution of the equation \eqref{generalformallydoublynonlinear} with $p$ and $m$ in the supercritical range \eqref{supercritical}, then $u$ is a locally bounded function, as shown in \cite{Vespri-Sosio-Fornaro2} and in \cite{Matias}.
\end{remark}

\section{ Short proof of Theorem \ref{T:1:1}} \label{S:3}

\subsubsection{The geometric setting}
Fix $(x_0,y_0) \in \Omega_T$ and construct the cylinder 
\begin{equation}
[(x_0,y_0)+Q(2\rho^p,2\rho^{\frac{p}{2}})] \subset \Omega_T,
\end{equation} \noindent After a translation we may assume that $(x_0,y_0)=(0,0)$. Let us set
$$\mu^+= \sup_{Q(\rho^p,\rho^{\frac{p}{2}})} u, \quad \mu^- = \inf_{Q(\rho^p,\rho^{\frac{p}{2}})} u, \quad \omega= \mu^+ - \mu^-$$
Consider the cylinder
\begin{equation} \label{1cylinder} Q(\rho^p, c_0 \rho), \quad \text{ where} \quad  c_0:= \omega^{\frac{p-2}{p}}
\end{equation} \noindent
To start the iteration, we assume that
\begin{equation}
\omega^{\frac{p-2}{p}} < \rho^{\frac{p-2}{2}}
\end{equation} \noindent 
otherwise if this is not the case, we would have
$$\omega \leq\rho^{\frac{p}{2}}.$$ \noindent
Thus we have
$$Q(\rho^p,c_0 \rho) \subset Q(\rho^p, \rho^{\frac{p}{2}}) \quad \text{and} \quad \essosc_{Q(\rho^p,c_0 \rho)} u \leq \omega $$ \noindent
Cylinders of the type of \eqref{1cylinder} have the space variables stretched by a factor $\omega $, which is intrinsically determined by the solution. If $p=2$ these are the standard parabolic cylinders.

\subsubsection{Transforming the variables and the PDE}
Introduce the change of variables
\begin{equation}
z=\frac{x}{c_o \rho }, \quad \tau= \frac{t}{\rho^{p}}, \quad v(z,\tau)= \frac{u(x,t)- \mu^-}{\omega }
\end{equation} \noindent
which maps
\begin{equation}
Q(\rho^p,c_0 \rho) \quad \rightarrow \quad B_1 \times (-1,0)
\end{equation} \noindent 
The transformed function $v$ solves an equation similar to \eqref{Eq:1:1}.
\subsubsection{Estimating positivity and conclusion}
Now we deal with the following two alternatives: either 
\begin{equation} \label{a}
|[v(z,0)>\frac{1}{2}] \cap B_{1/2}| > \frac{1}{2} |B_{1/2}|
\end{equation} \noindent
or we would have 
\begin{equation} \label{b}
|[v(z,0)>\frac{1}{2}] \cap B_{1/2}| \leq \frac{1}{2} |B_{1/2}|
\end{equation} \noindent and
as the function $(1-v)$ still satisfies equation \eqref{Eq:1:1} with similar structure conditions, we can assume that \eqref{a} holds. Thus we suppose \eqref{a} and by Propostion 2.1 we have for $t_o <0$
\begin{equation} \label{conto1}\frac{1}{2^{N+2}} |B_1| = \frac{1}{4} |B_{\frac{1}{2}}|  \leq  \sup_{t_0 \leq \tau \leq 0}  \bigg( \int_{[v>\frac{1}{2}] \cap B_{1/2}} v(z, \tau) dz + \int_{[v \leq 1/2] \cap B_{1/2}} v(z,\tau) dz   \bigg) $$$$ 
\leq \gm \bigg{\{} \inf_{t_0 \leq \tau \leq 0} \int_{B_1} v(z,\tau) dz + \bigg( \frac{|t_o|}{( 1/2)^{N(p-2)+p}} \bigg)^{\frac{1}{2-p}} \bigg{\}}.
\end{equation}
If we take 
\begin{equation}  |t_0|^{\frac{1}{2-p}} \leq \frac{1}{\gm 2^{N+3}} 2^{(N-\frac{p}{2-p})} |B_1| = \frac{1}{\gm 2^{3+\frac{p}{2-p}}} |B_1| \end{equation} which can be done by defining 
\begin{equation} \label{t0} t_0 = -\bigg( \frac{1}{\gm 2^{3+\frac{p}{2-p}}}\bigg)^{2-p} \end{equation}
we obtain the information
\begin{equation} \label{conto2}
\inf_{t_0 \leq \tau \leq 0} \int_{B_1} u(x,\tau) \, dx \ge \frac{1}{ \gm 2^{N+3}} |B_1| = \frac{2^{\frac{p}{2-p}}}{2^N} |t_0|^{\frac{1}{2-p}} |B_1| = 4 \eta \, |B_1|
\end{equation} \noindent
Where we have defined
$$\eta=\frac{2^{\frac{p}{2-p}}}{2^{N+2}} |t_0|^{\frac{1}{2-p}}$$
This implies that
\begin{equation}
    |[u>\eta] \cap B_1|> \eta \, |B_1|, \quad \text{for all} \quad  \tau \in (t_0,0]
    \end{equation} \noindent
We apply Proposition 1.1 with 
$$s=0,\quad  M= \eta, \quad \quad \epsilon= \frac{t_0}{\eta^{2-p}}= \frac{2^p}{2^{\frac{N+2}{2-p}}},\quad \quad B_{16} \times (t_0,0] \subset D_{v}  $$ being $D_v$ the domain of $v$ function, to get that there exists a $\sigma \in (0,1) $ and $ \epsilon^* \in (0, \frac{\epsilon}{2}]$\\
such that

\begin{equation}
v(z,\tau)  \ge \sigma \eta , \quad \text{for all} \quad z \in B_{2} 
\end{equation} \noindent for all times
\begin{equation}
    - \epsilon_1 t_0 = - \epsilon^* \frac{2^p}{2^{\frac{N+2}{2-p}}} \,
    t_0 \leq \tau \leq 0
\end{equation} \noindent 
Returning back to the original coordinates this means that
\begin{equation}
u(x,t) \ge \mu^{-} + \sigma \eta\,  \omega , \quad \forall x \in B_{c_o \rho}, \quad \sigma, \eta \in (0,1) 
\end{equation} \noindent for all times
\begin{equation}
- \epsilon_1 \,  t_0 \rho^p \leq t \leq 0 
\end{equation} \noindent
This implies
\begin{equation}
\essosc_{Q((\frac{\rho}{2^l})^p, c_0 \rho)} u \leq (1- \sigma \eta ) \omega
\end{equation}
for $l= \frac{1}{p} \log_2 \bigg( \frac{1}{\epsilon_1 t_0} \bigg)$ given by the request $2^{-lp}= \epsilon_1 t_0 $. We are in the hypothesis of Proposition \ref{P:iteration}, as the process can now be repeated inductively starting from such relation.

\section{Short proof of Theorem \ref{HCDBNL}} \label{S:4}

\subsection{Geometrical setting and the alternative}
Define $$M = \sup_{Q_{\rho}(\theta)} u $$
and begin by normalizing the function by the transformation
\begin{equation} \label{transf}
v(x,t)= \frac{u(x,t)}{M}, \quad \quad 0 \leq v \leq 1
\end{equation} \noindent
Let $0<\epsilon_0<1$ be a number to be defined later in \eqref{epsilonzero}, and consider the following cases:
if 
\begin{equation} \label{inf1}
\inf_{Q_{\rho}(\theta)} v \ge \epsilon_0
\end{equation} \noindent
then the equation \eqref{doublynonlinear} behaves as a variable coefficients $p$-laplacean type equation, and by arguments of previous \textsection \ref{S:3} we have the reduction of oscillation. If otherwise
\begin{equation} \label{inf2}
\inf_{Q_{\rho}(\theta)} v < \epsilon_0
\end{equation} \noindent we may suppose the worst case, which is 
\begin{equation*}
\inf_{Q_{\rho}(\theta)} v =0
\end{equation*} \noindent 
Finally we set the alternative on the measure of the positivity set of $v$. We set $M=\frac{1}{2}$, and consequently $\theta=1$. Let us suppose that for a sufficiently small $\rho$ to be fixed later, that
$$Q_1 =B_{1} \times (- 1,0] \subset B_{4} \times (- 4^p,0] \subset  \Omega_T$$
If $\nu \in (0,1)$ is the number of Lemma \ref{DG}, we can set two alternatives: either
\begin{equation} \label{alt1}
|[v(x,t)\ge \frac{1}{2}] \cap Q_{\rho}(\theta)| \ge \nu |Q_{\rho}(\theta)| 
\end{equation} \noindent or
\begin{equation} \label{alt2bis}
|[v(x,t) > \frac{1}{2}] \cap Q_{\rho}(\theta)| < \nu | Q_{\rho}(\theta)|
\end{equation} \noindent
\subsection{Conclusion of the proof of the Theorem \ref{HCDBNL}}
Assume \eqref{alt1} holds, then we have that it exists a $ \bar{t} \in (- \rho^p,0]$ such that
\begin{equation}
|[v(x,\bar{t})>\frac{1}{2}] \cap B_{1}|\ge |[v(x,\bar{t})>\frac{1}{2}] \cap B_{\rho}| > \nu \rho^N |B_{1}|= \nu w_N \rho^N
\end{equation}
By $L^1$-Harnack inequality applied in the box $B_1$ and by estimating $T \leq \rho^p$ we have that 
\begin{equation*}
 \frac{\nu w_N}{2} \rho^N \leq \int_{B_1} v(x,\bar{t}) dx \leq    \sup_{-\rho^p \leq \tau \leq 0} \int_{B_{2}} v(x,\tau) dx \leq \gm \bigg{\{} \inf_{\rho^p \leq \tau \leq  0} \int_{B_2} v(x,\tau) dx + ( \rho^p )^{\frac{1}{3-m-p}}  \bigg{\}}
\end{equation*} \noindent
So, by asking the condition of supercritical range $m+p>3-\frac{p}{N}$ we have
\begin{equation*} \label{condition}
\gm (\rho^p)^{\frac{1}{3-m-p}} \leq \frac{ \nu w_N \rho^N }{4}, \quad  \quad \rho \leq \bigg( \frac{\nu  w_N}{4 \gm} \bigg)^{\frac{3-m-p}{p-N(3-m-p)}} =:\rho_0
\end{equation*} \noindent denoting with $w_N$ the Lebesgue measure of the $N$-dimensional ball, we have
\begin{equation} 
\inf_{-\rho_0^p \leq \tau \leq 0} \int_{B_{2}} v(x,\tau) dx \ge \frac{\nu w_N}{\gm 4} \rho_0^N= \bigg( \frac{\nu w_N}{4 \gm}\bigg)^{\frac{p}{p-N(3-m-p)}}=: \eta_1
\end{equation} \noindent 
This implies 
\begin{equation*}
|[v(x,t)>\frac{\eta_1}{2}] \cap B_2| > \frac{\eta_1}{2} |B_{2}|, \quad \text{for all} \quad t \in (-\rho_0^p,0  ]
\end{equation*} \noindent
Finally we use expansion of positivity Lemma to get 
\begin{equation} \label{positivity-alt1}
v(x,t) \ge \eta \eta_1= \eta \bigg( \frac{\nu w_N}{4 \gm}\bigg)^{\frac{p}{p-N(3-m-p)}}=: \eta^* , \quad \text{in} \quad B_4 \times (- \rho_0^p,0 ]
\end{equation} \noindent
Now we can choose 
\begin{equation} \label{epsilonzero}
\epsilon_0=  \frac{\eta^*}{2}
\end{equation} \noindent 
where $\epsilon_0$ is the constant defined in \eqref{inf1}.$$$$
If otherwise \eqref{alt2bis} holds, we use Lemma \ref{DG} for which we take 
$$M=\frac{1}{2} , \quad \quad \theta=(2M)^{3-m-p}, \quad \quad c_0=1$$ we fulfill its hypothesis to have
\begin{equation}
v(x,t) \leq \frac{1}{2} \sqrt[\beta]{\frac{3}{2}} \quad \text{in} \quad B_{\rho_0 /2} \times \bigg( \bigg(-\frac{\rho_0}{2}\bigg)^p \bigg( \frac{1}{2}\bigg)^{3-m-p},0 \bigg]=Q_{\frac{\rho_0}{2}}(\theta)
\end{equation} \noindent 
$$$$
Finally, if \eqref{inf2} holds, then by expansion of positivity we have demonstrated that in the two alternatives \eqref{alt1} and \eqref{alt2bis} we obtain respectively 
\[
\inf_{Q_{\frac{\rho_0}{2}}} v \ge 2 \epsilon_0 \quad \text{or} \quad \sup_{Q_{\frac{\rho_0}{2}}}v \leq \frac{1}{2} \sqrt[\beta]{\frac{3}{2}} 
\] 
while, if \eqref{inf1} holds we have an equation of the $p$-Laplacean type and by the same technique of previous section we still arrive to an estimate of the previous kind.
In either case, returning to the original function, we obtain a reduction of oscillation and therefore the H\"older continuity in a similar fashion than we did in the previous section (we refer for details to \cite{Matias}). 

\subsection*{Acknowledgements:} 
We wish to thank Matias Vestberg and Naian Liao for helpful conversations on the subject. Moreover, both authors are partially founded by INdAM (GNAMPA).

\end{document}